\input amstex
\input Amstex-document.sty
\loadeufm
\loadmsbm
\loadeufm

\catcode `\@=11
\def\logo@{}
\catcode `\@=12
%\magnification \magstep1
%\NoRunningHeads
%\NoBlackBoxes
%\TagsOnLeft

\def \={\ = \ }
\def \+{\ +\ }
\def \-{\ - \ }

\def \b|{\big |}

\def \g1{\Gamma_1}

\def\rarr#1#2{\smash{\mathop{\hbox to .5in{\rightarrowfill}}
     \limits^{\scriptstyle#1}_{\scriptstyle#2}}}

\def\larr#1#2{\smash{\mathop{\hbox to .5in{\leftarrowfill}}
      \limits^{\scriptstyle#1}_{\scriptstyle#2}}}

\def\swarr#1#2 {\llap{$\scriptstyle #1$}  \swarrow
    \vcenter to .5in{}\rlap{$\scriptstyle #2$}}

\pageno 87

\def\shorttitle{Diophantine Geometry over Groups  and the
Elementary Theory $\cdots$}
\def\shortauthor{Z. Sela}

\topmatter
\title\nofrills{\boldHuge Diophantine Geometry over Groups  and the
Elementary Theory of Free and Hyperbolic Groups*}
\endtitle

\author \Large Z. Sela$^\dag$
\endauthor

\thanks *Partially supported by an Israel academy of sciences fellowship, an NSF grant DMS9729992
through the IAS, and the IHES. \endthanks \abstract\nofrills \centerline{\boldnormal Abstract}

\thanks $\dag$Hebrew University, Jerusalem 91904, Israel. E-mail: zlil\@math.huji.ac.il
\endthanks

\vskip 4.5mm

{\ninepoint We study sets of solutions to equations over a free group, projections of such sets, and the structure
of elementary sets defined over a free group. The structre theory we obtain enable us to answer some questions of
A. Tarski's, and classify those finitely generated groups that are elementary equivalent to a free group.
Connections with low dimensional topology, and a generalization to (Gromov) hyperbolic groups will also be
discussed.

\vskip 4.5mm

\noindent {\bf 2000 Mathematics Subject Classification:} 14, 20.}

%\noindent {\bf Keywords and Phrases:} Cohomology, Symmetric group, Rotation group.}
\endabstract
\endtopmatter

\document

\baselineskip 4.5mm \parindent 8mm

Sets of solutions to equations defined over a free group have been studied
extensively, mostly since Alfred Tarski presented his fundamental
questions on the
elementary theory of free groups in the mid 1940's. Considerable progress in the study of
such sets of solutions was made by G. S. Makanin, who constructed
an algorithm that decides if a system of equations defined over a free group
has a solution [Ma1], and showed that the universal and positive theories
of a free group are decidable [Ma2]. A. A. Razborov
was able to give a description of the entire set of solutions to a system of
equations defined over a free group [Ra], a description that was further developed
by O. Kharlampovich and A. Myasnikov [Kh-My].

A set of solutions to equations defined over a free group
is clearly a discrete set, and all the previous techniques and methods
that studied these sets
are combinatorial in nature.
Naturally, the structure of sets of solutions defined over a free group is
very different from the structure of sets of solutions (varieties)
to systems of equations  defined over the complexes, reals or a number field.
Still, perhaps surprisingly, concepts from complex algebraic geometry and
from Diophantine geometry can be borrowed to study varieties defined over
a free group.

In this work we borrow concepts and techniques
from geometric
group theory, low dimensional topology, and Diophantine geometry to study
the structure of varieties defined over a free (and hyperbolic) group. Our techniques and
point of view on the study of these varieties is rather different from any
of the pre-existing techniques in this field, though, as one can expect,
some of our preliminary results overlap with previously known ones.
The techniques and concepts we use
enable the study of the structure of varieties defined over a free group and their
projections (Diophantine sets), and in particular, give us
the possibility to answer some questions that  seem to be essential in
any attempt to understand the structure of elementary sentences and predicates
defined over a free (and hyperbolic) group.

In this note we summarize the main results of our work, that enable one to answer affirmatively some of A.
Tarski's problems on the elementary theory of a free group, and  classify those finitely generated groups that are
elementary equivalent to a (non-abelian) free group. we further survey some of our results on the elementary
theory of a (torsion-free) hyperbolic group, that generalize the results on free groups. The work itself appears
in [Se1]-[Se8].

We start with what we see as the main result on the elementary theory
of a free group we obtained - quantifier elimination. Quantifier elimination and its proof
is behind all the other results presented in this note.

\proclaim{Theorem 1 ([Se7],1)} Let $F$ be a non-abelian free group, and let
$Q(p)$ be a definable set over $F$. Then $Q(p)$ is in the Boolean
algebra of $AE$ sets over $F$.
\endproclaim

In fact it is possible to give a strengthening of
theorem 1 that specifies a subclass of $AE$ sets that generates the Boolean algebra
of definable sets, a more refined description that is essential in studying other model-theoretic
properties of the elementary theory of a free group.

Theorem 1 proves that every definable set over a free group is in the Boolean algebra of $AE$ sets. To answer
Tarski's questions on the elementary theory of a free group, i.e., to show the equivalence of the elementary
theories of free groups of various ranks, we need to show that for coefficient free predicates, our quantifier
elimination procedure does not depend on the rank of the coefficient group.

\proclaim{Theorem 2 ([Se7],2)} Let $Q(p)$ be a set defined by a coefficient-free predicate over a
group. Then there exists a set $L(p)$ defined by a coefficient-free predicate which is in the Boolean algebra of
$AE$ predicates, so that for every non-abelian free group $F$, the sets $Q(p)$ and $L(p)$ are equivalent.
\endproclaim

Theorem 2 proves that in handling coefficient-free predicates,
our quantifier elimination procedure does not depend
on the rank of the coefficient (free) group. This together with the equivalence of the
$AE$ theories of free groups ([Sa],[Hr]) implies
an affirmative answer to Tarski's problem on the equivalence of the
elementary theories of  free groups.

\proclaim{Theorem 3 ([Se7],3)}  The elementary theories of non-abelian  free groups
are equivalent.
\endproclaim

Arguments similar to the ones used to prove theorems 2 and 3, enable us
to answer affirmatively another question of Tarski's.

\proclaim{Theorem 4 ([Se7],4)} Let $F_k,F_{\ell}$ be free groups for
$2 \leq k \leq \ell$. Then the standard embedding $F_k \to F_{\ell}$
is an elementary embedding.

More generally, let $F,F_1$ be non-abelian free groups, let $F_2$ be a free group, and suppose that $F=F_1*F_2$.
Then the standard embedding $F_1 \to F$ is an elementary embedding.
\endproclaim

Tarski's problems deal with the equivalence of the elementary theories of free
groups of different ranks. Our next goal is to
get a classification of all the f.g.\ groups that
are elementary equivalent to a free group.

Non-abelian $\omega$-residually free groups (limit groups) are known to be the f.g.\ groups
that are universally equivalent to a
non-abelian free group. If a limit group contains a free abelian group of rank 2,
it can not be elementary equivalent to a free group. Hence, a f.g.\ group that
is elementary equivalent to a non-abelian free  group must be a non-elementary
(Gromov) hyperbolic  limit group. However, not every non-elementary hyperbolic limit group
is elementary equivalent to a free group.

\noindent
To demonstrate that we look at the following example.
Suppose that $G=F*_{<w>} \, F=<b_1,b_2>*_{<w>} \, <b_3,b_4>$ is a double of a free group of rank 2,
suppose that $w$ has no roots in $F$, and suppose that the given amalgamated product is the
abelian JSJ decomposition of the group $G$. By our assumptions, $G$ is  a hyperbolic limit group (see [Se1], theorem 5.12).

\proclaim{Claim 5 ([Se7],5)} The group $G=F*_{<w>} \, F$ is not elementary equivalent to the free group $F$.
\endproclaim

In section 6 of
[Se1] we have presented $\omega$-residually free towers, as an example of limit groups  (the same
groups  are  presented in [Kh-My] as well, and are called there
NTQ groups).

A hyperbolic $\omega$-residually free tower is constructed in finitely many steps. In its first level there is  a
non-cyclic free product of (possibly none) (closed) surface groups and a (possibly trivial) free group, where each
surface in this free product is a hyperbolic surface (i.e., with negative Euler characteristic), except the
non-orientable surface of genus 2. In each additional level we add a punctured surface that is amalgamated to the
group associated with the previous levels along its boundary components, and in addition there exists a retract
map of the obtained group onto the group associated with the previous levels. The punctured surfaces are supposed
to be of Euler characteristic bounded above by -2, or a punctured torus.

The procedure used for eliminating quantifiers over a free group enables us to show that
every hyperbolic $\omega$-residually free tower is elementary equivalent to a free group.
The converse is obtained by using basic properties of the JSJ decomposition and
the (canonical) Makanin-Razborov diagram of a limit group ([Se7], theorem 6).
Therefore, we are finally able to get
a classification of those f.g.\ groups that are elementary equivalent to a free group.

\proclaim{Theorem 6 ([Se7],7)} A f.g.\ group is elementary equivalent to a
non-abelian  free group if and only if it is a non-elementary
hyperbolic $\omega$-residually free tower.
\endproclaim

So far we summarized the main results of our
work, that enable one to answer affirmatively
some of A. Tarski's problems on the elementary theory
of a free group, and  classify those finitely generated groups that are elementary
equivalent to a (non-abelian) free group. In the rest of this note
we survey some of our results on the elementary theory of a (torsion-free) hyperbolic
group, that generalize the results presented  for a free group.

In the case of a free group, we have shown that every definable set is in the Boolean
algebra of $AE$ sets. The same holds for a general hyperbolic group.

\proclaim{Theorem 7 ([Se8],6.5)} Let $\Gamma$ be a non-elementary torsion-free
hyperbolic group, and let
$Q(p)$ be a definable set over $\Gamma$. Then $Q(p)$ is in the Boolean
algebra of $AE$ sets over $\Gamma$.

Furthermore, if $Q(p)$ is a set defined by a coefficient-free predicate defined over  $\Gamma$, then $Q(p)$ can be
defined by a coefficient-free predicate which is in the Boolean algebra of $AE$ predicates.
\endproclaim

The procedure used for quantifier elimination over a free group
enabled us to get a classification of
those f.g.\ groups that are elementary equivalent to a free group (theorem 6). In a similar
way, it is possible to get a classification of those f.g.\ groups that are elementary
equivalent to a given torsion-free hyperbolic group.

\noindent
We start with the following basic fact, that shows  the elementary invariance of
 negative curvature in groups.

\proclaim{Theorem 8 ([Se8],7.10)} Let $\Gamma$ be a torsion-free hyperbolic group, and
let $G$ be a f.g.\ group. If $G$ is elementary equivalent to $\Gamma$, then
$G$ is a torsion-free hyperbolic group.
\endproclaim

Theorem 8 restricts the class of f.g.\ groups that are elementary equivalent to a given
hyperbolic group, to the class of hyperbolic groups. To present the elementary classification of
hyperbolic groups we start with the following basic fact.

\proclaim{Proposition 9 ([Se8],7.1)} Let $\Gamma_1,\Gamma_2$ be non-elementary torsion-free
rigid
hyperbolic groups  (i.e., $\Gamma_1$ and $\Gamma_2$ are freely-indecomposable
and do not admit any non-trivial cyclic splitting).
Then $\Gamma_1$ is elementary equivalent to $\Gamma_2$ if and only if
$\Gamma_1$ is
isomorphic to $\Gamma_2$.
\endproclaim

Proposition 9 implies that, in particular, a uniform lattice  in a real rank
1 semi-simple Lie group that is not $SL_2(R)$ is elementary equivalent to another such lattice
if and only if the two lattices are isomorphic, hence, by Mostow's rigidity
the two lattices are conjugate in the same Lie group.
By Margulis's
normality and super-rigidity theorems, the same hold in higher rank (real)
Lie groups.

\proclaim{Theorem 10 ([Se8],7.2)} Let $L_1,L_2$ be uniform lattices in real semi-simple Lie
groups that are not $SL_2(R)$.
Then $L_1$ is elementary equivalent to $L_2$ if and only if
$L_1$ and $L_2$ are conjugate lattices in the same real Lie group $G$.
\endproclaim

Proposition 9 shows that rigid hyperbolic groups are elementary equivalent
if and only if they are isomorphic. To classify elementary equivalence
classes of hyperbolic groups in general, we associate with every (torsion-free)
hyperbolic group $\Gamma$, a subgroup of it, that we call the $elementary$
$core$ of $\Gamma$, and denote $EC(\Gamma)$. The elementary core is a retract of
the ambient hyperbolic group $\Gamma$, and although it is not canonical, its isomorphism
type is an invariant of the ambient hyperbolic group. The elementary core
 is constructed iteratively from the ambient hyperbolic group as we describe in definition
7.5 in [Se8].

\noindent
The elementary core of a hyperbolic group is a prototype for its elementary theory.

\proclaim{Theorem 11 ([Se8],7.6)} Let $\Gamma$ be a non-elementary torsion-free hyperbolic group that is not a $\omega$-residually free tower, i.e., that is not
elementary equivalent to a free group. Then $\Gamma$ is elementary equivalent to its elementary core $EC(\Gamma)$. Furthermore,
the embedding of the elementary core $EC(\Gamma)$ in the ambient group $\Gamma$ is an elementary
embedding.
\endproclaim

Finally, the elementary core is a complete invariant of the class of groups that are elementary
equivalent to a given (torsion-free) hyperbolic group.

\proclaim{Theorem 12 ([Se8],7.9)} Let $\Gamma_1,\Gamma_2$ be two  non-elementary torsion-free hyperbolic
groups.
Then $\Gamma_1$ and $\Gamma_2$ are elementary equivalent if and only if their elementary cores $EC(\Gamma_1)$ and  $EC(\Gamma_2)$ are isomorphic.
\endproclaim

Theorem 12 asserts that the elementary class of a torsion-free hyperbolic group is determined by the isomorphism
class of its elementary core. Hence, in order to be able to decide whether two torsion-free hyperbolic groups are
elementary equivalent one needs to
  compute their elementary core, and to decide if the two elementary
cores are
isomorphic. Both can be done using the solution to the isomorphism problem for torsion-free
hyperbolic groups.

\proclaim{Theorem 13 ([Se8],7.11)} Let $\Gamma_1,\Gamma2$ be two  torsion-free hyperbolic groups. Then it is decidable if $\Gamma_1$ is elementary equivalent to $\Gamma_2$.
\endproclaim

%\newpage

\specialhead \noindent \boldLARGE References \endspecialhead

\widestnumber\key{XX-XXX}

\ref\key Hr
\by E. Hrushovski
\paper private communication
\endref

\ref\key  Kh-My \by O. Kharlampovich and A. Myasnikov \paper
Irreducible affine varieties over a free group II \jour Jour. of
Algebra \vol 200 \yr 1998 \pages 517--570
\endref

\ref\key  Ma1 \by G. S. Makanin \paper Equations in a free group
\jour Math. USSR Izvestiya \vol 21 \yr 1983 \pages 449--469
\endref

\ref\key  Ma2 \bysame \paper Decidability of the universal and
positive theories of a free group \jour Math. USSR Izvestiya \vol
25 \yr 1985 \pages 75--88
\endref

\ref\key Ra
\by A. A. Razborov
\paper On systems of equations in a free group
\paperinfo Ph.D. thesis, Steklov Math. institute, 1987
\endref

\ref\key Sa \by G. S. Sacerdote \paper Elementary properties of
free groups \jour Transactions Amer. Math. Soc. \vol 178   \yr
1973   \pages 127--138
\endref

\ref\key Se1-Se8 \by Z. Sela \paper Diophantine geometry over groups I-VIII \paperinfo preprints, \linebreak
www.ma.huji.ac.il/~zlil
\endref

\enddocument